\documentclass[11pt]{amsart}
\usepackage{amsthm}
\usepackage{graphicx}
\usepackage{bm}
\usepackage{graphicx}
\usepackage{amssymb}
\usepackage{amsmath}
\usepackage{enumerate}
\newtheorem{thm}{Theorem}
\newtheorem{lem}{Lemma}

\title{Inverse Nodal Problems}
\author{David Klawonn}
\address{David.Klawonn@aei.mpg.de, Max Planck Institute for Gravitational Physics (Albert Einstein Institute), Am Muehlenberg 1, D-14476 Golm, Germany}

\begin{document}
\maketitle

\begin{abstract}
It is shown that nodal sequences determine the underlying manifold up to scaling within classes of rectangles with Dirichlet boundary conditions, separable two dimensional tori, two-dimensional flat Klein bottles and flat tori in two and three dimensions. 
\end{abstract}

\section{Introduction}
One of the most famous \emph{Inverse Problems} in mathematical physics is the question whether one can hear the shape of a drum, which was posed by \emph{Kac} in 1966 \cite{kac}. The spectrum of the Laplace operator $\sigma_{M}(\triangle)$ on a Riemannian manifold $M$, representing the drum, corresponds to the frequencies of the eigenmodes in which this manifold is vibrating. Since 1964 it is known that the answer to this question is negative, when \emph{Milnor} \cite{milnor} found a pair of isospectral, non-isometric flat tori in 16 dimensions. Many examples of such pairs are known today. For an overview on different aspects of the theory of isospectral manifolds the works of \emph{Sunada}, \emph{Gordon}, \emph{Webb}, \emph{Wolpert} and \emph{Brooks} should be mentioned \cite{SU85,GOWEWO92,BR99}. Recent results and examples of isospectral manifolds are presented in \cite{BAPA08,BAPABE08}. Forty years after the foundation of the area of \emph{Isospectrality} another invariant object of the manifold $M$, connected to its vibration modes, has joined the field of interest \cite{blumetal02}: The eigenfunctions corresponding to a specific eigenvalue $\lambda$ form an \emph{eigenspace} $\mathfrak{E}_\lambda$, due to the equation $\triangle \Psi = \lambda\Psi$. For every element of the eigenspace $\Psi\in\mathfrak{E}_\lambda$ its \emph{nodal domains} are defined as the connected components of the set $M\setminus \Psi^{-1}(0)$. The set $\Psi^{-1}(0)$ is called the \emph{nodal set of $\Psi$}. A visualization of the nodal set is given by the \emph{Chladni figures}, that arise as those patterns formed by sand poured on a vibrating membrane \cite{ullmann83}.
Assigning to each eigenvalue $\lambda$ the set of \emph{nodal counts} $\nu(\Psi)$, the numbers of nodal domains $\mathcal{N}(\lambda)=\{\nu(\Psi_1),\nu(\Psi_2),\dots\}\footnote{The definition of $\mathcal{N}(\lambda)$ is specific to each class of manifolds. This will be done later in this note.}$, that are obtained by elements of its eigenspace $\Psi_1,\Psi_2,\dots\in\mathfrak{E}_\lambda$ and ordering these sets according to the position of $\lambda$ within the spectrum, we obtain a sequence
$
\mathfrak{N}(M)=\{\mathcal{N}(\lambda_1),\,\mathcal{N}(\lambda_2),\,\hdots\, ,\mathcal{N}(\lambda_k),\,\hdots\}
$
called the \emph{nodal sequence}.\\
\indent
Due to the existence of isospectral manifolds, it is known that the information stored in the spectrum of a Riemannian manifold does not suffice to determine the manifold itself. It is therefore interesting to investigate whether the nodal sequence reveals similar or complementary information about the underlying geometry. This idea is due to \emph{Smilansky} (cf. e.g. \cite{blumetal02, gnutzkara06}).  Recent progress shows, that this approach is very promising. For different objects like graphs, billiards, tori and surfaces of revolution, nodal sequences have been investigated and in each case they reveal information about the underlying objects. It was possible to distinguish between integrable and chaotic billiards by using the nodal sequence \cite{blumetal02}. Examples of isospectral manifolds and graphs have been distinguished by using their nodal sequences \cite{BASHSM06,SmilGnutz,BRKLPU07}, which leads to the assumption that the nodal sequence contains complementary information about the geometry to the spectrum. It was even possible to find a trace formula for the nodal counting function of different manifolds \cite{gnutzkara06}, which shows a deep similarity to the ideas that once led to the original isospectrality problem.\\
\indent
We want to focus now on a specific part of the investigation of spectra of manifolds. Changing the isospectrality problem slightly by restricting it to a class of manifolds leads to a new question. Are there classes within which it is possible to hear the shape of a drum? The answer is positive, and examples for such classes are low-dimensional flat tori or surfaces of revolution \cite{BRHE84,ZE98,BEGAMA71}.
This problem can be translated to the study of nodal sequences, which we call the \emph{Inverse Nodal Problem}.
It was solved recently by \emph{Smilansky et al.} for surfaces of revolution \cite{karasmil08}, and rectangles were treated in \cite{SMSA05}.\\
\indent
In this note methods to reconstruct several types of manifolds from their nodal sequences are presented.  It is always assumed that the given nodal sequence belongs to a manifold of a specific type. We will investigate two- and three-dimensional flat tori. Furthermore we will study two dimensional rectangles, separable tori and Klein bottles. Within all these classes the \emph{Inverse Nodal Problem} will be solved.\\
\\
\indent The work supports the conjecture indicated by \emph{Smilansky et al.} about the importance of the nodal sequence for the determination of manifolds by their vibration behaviour, and therefore should encourage the further study of nodal sequences in general. 


\section{Some Manifolds and their Nodal Sequences}
In this section different classes of manifolds are introduced. Formulas for the number of nodal domains are given. As the number of nodal domains is invariant under scaling of the manifold, we will reduce the free parameters in each class of manifolds by one. All manifolds are equipped with the euclidean metric.
\subsection{Dirichlet-rectangles and 2-dimensional separable tori}
We define $R_{\alpha}=[0,a] \times [0,b]$  for $a,b \in \mathbb{R_+}$ and set $\alpha b = a=1$ to reduce the free parameters by one.
Without loss of generality we choose $\alpha\leq1$, since $R_{\alpha}$ and $R_{\alpha^{-1}}$ are isometric. The corresponding class is called $\mathfrak{R}=\{R_\alpha\,:\,\alpha\leq1\}$. Considering the standard Laplacian $\triangle$ on $R_\alpha\in\mathfrak{R}$, the Dirichlet problem on $R_{\alpha}$ is given by
\[
\triangle \varphi = \lambda \varphi, \quad \varphi |_{\partial R_{\alpha}}=0,\,\varphi\in C^{\infty}(R_\alpha),\, \lambda\in\mathbb{C}\,.
\]
With $m,n\in\mathbb{N}_0$ the solutions are 
\[
\varphi_{mn}=\sin\left(\pi m  x/a\right)\sin\left(\pi n  y/b\right) \mbox{ with } \lambda_{mn}=\pi^2\left(m^2/a^2+n^2/b^2\right)=\left(\pi/a\right)^2\left(m^2+\alpha^2n^2\right)\,.
\] 
As we set $a=1$ we obtain
\[
\lambda_{mn}=\pi^2\left(m^2+\alpha^2n^2\right)\,.
\]
The nodal lines of the eigenfunctions form a regular checkerboard pattern and the number of nodal domains is therefore given by
\begin{equation}\label{nod-1}
\nu_{mn}=\nu(\varphi_{mn})=m\cdot n\,.
\end{equation}
To introduce the class of two-dimensional separable flat tori, for $a\leq b\in\mathbb{R}_+$ with $\alpha b=a=1$, the equivalence relation $\sim_\alpha^T$ in $\mathbb{R}^2$ is defined:
\[
\mathbb{R}^2\ni(x,y)\sim^T_{\alpha}(x',y') \Longleftrightarrow (x',y')=(x+k\cdot a,y+l\cdot b) \mbox{ for } a,\,b\in\mathbb{R}_+\mbox{ and }k,l\in\mathbb{Z}\,.
\]
We define the corresponding flat torus $T_\alpha$ by $T_{\alpha}=\mathbb{R}^2/\sim_{\alpha}^T$.
The corresponding class is called $\mathfrak{T}=\{T_{\alpha}\,:\,\alpha\leq 1\}$.
Considering the standard Laplacian $\triangle$ on $T_{\alpha}$, the eigenfunctions are the solutions of 
\[
\triangle\varphi=\lambda\varphi,\quad \varphi(x,y)=\varphi(x+k\cdot a,y)=\varphi(x,y+k\cdot b),\,\varphi\in C^{\infty}(\mathbb{R}^2),\,\lambda\in\mathbb{C},\, k,l\in\mathbb{Z}\,.
\]
For $m,n\in\mathbb{Z}$ we obtain the following solutions:
\[
\varphi_{mn}(x,y)=\begin{cases}
\cos(2\pi mx/a)\cos(2\pi ny/b)&m,n\geq0\\
\sin(2\pi mx/a)\sin(2\pi ny/b)&m,n<0\\
\cos(2\pi mx/a)\sin(2\pi ny/b)&m\geq0,\,n<0\\
\sin(2\pi mx/a)\cos(2\pi ny/b)&m<0,\,n\geq0\,.
\end{cases}
\]
The eigenvalues are given by
\begin{eqnarray*}
\lambda_{mn}&=&4\left(\pi / a\right)^2\left(m^2+\alpha^2n^2\right)\\
&=&4\pi^2\left(m^2+\alpha^2n^2\right)
\end{eqnarray*}
The structure of the nodal set is similar to the rectangle case. According to \cite{gnutzkara06} the number of nodal domains is given by
\begin{equation}\label{nod-2}
\nu_{mn}=(2|m|+\delta_{m0})(2|n|+\delta_{n0})\,.
\end{equation}

\subsection{Flat Klein bottles}
\begin{figure}[htbp]
\begin{center}
\includegraphics[width=8cm]{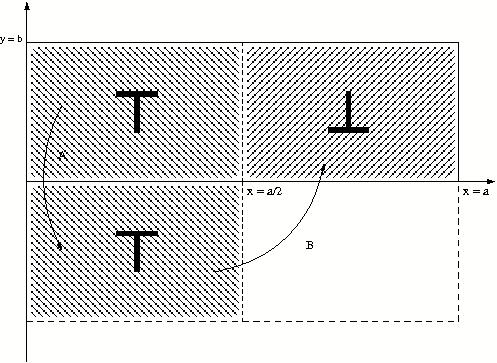}\\ 
\caption{{\bf Klein bottle}}
\label{KB}
\end{center}
\end{figure}

Given $a,b\in\mathbb{R}_+$ with $\alpha b = a=1$ and the following equivalence condition
\begin{eqnarray*}
(x,y)\sim^K_\alpha (x',y') &\Longleftrightarrow& (x'=x \mbox{ and } y'=y+k\cdot b)\qquad \mbox{ or } \\
 && (x'=x+l\cdot a/2 \mbox{ and } y'=-y)\,\,\mbox{ or }\\
&& (x'=x+m\cdot a\mbox{ and }y'=y)\,\,\,\mbox{ for }k\in\mathbb{Z},\,l\in2\mathbb{Z}+1,\,m\in\mathbb{Z} \,.
\end{eqnarray*}
The Klein bottle $K_{\alpha}$ is defined by $K_\alpha=\mathbb{R}^2/\sim_{\alpha}^K$. The corresponding class is given by $\mathfrak{K}=\{K_\alpha\,:\,\alpha\in\mathbb{R}_+\}$. A visualization of the Klein bottle is given in figure \ref{KB}, where the symmetry operations are indicated by $A$ and $B$.
Considering the standard Laplacian $\triangle$, the eigenfunctions are the solutions of the following problem:
\begin{eqnarray*}
&&\triangle\varphi=\lambda\varphi,\, \varphi(x,y)=\varphi(x,y+k\cdot b)=\varphi(x+l\cdot a/2,-y)=\varphi(x+m\cdot a,y),\\
&&\mbox{where }\varphi\in C^{\infty}(\mathbb{R}^2),\,\lambda\in\mathbb{C},\,k,m\in\mathbb{Z}, l\in2\mathbb{Z}+1\,.
\end{eqnarray*}

The corresponding eigenfunctions for $\lambda_{mn}$ are then given by
\[
\Psi_{mn}(x)=\begin{cases}
\exp(2\pi i mx/a)\cos(2\pi ny/b)\quad m \mbox{ even, } n\geq 0\\
\exp(2\pi i  mx/a)\sin(2\pi ny/b)\quad m \mbox{ odd, } n\neq 0\,.
\end{cases}
\]
Considering only the real parts we obtain eigenfunctions of the form
\[
\Psi^{\mathfrak{Re}}_{mn}(x)=\begin{cases}
\cos(2\pi mx/a)\cos(2\pi ny/b)\quad m \mbox{ even, } n\geq 0\\
\cos(2\pi  mx/a)\sin(2\pi ny/b)\quad m \mbox{ odd, } n\neq 0\,.
\end{cases}
\]

The eigenvalues are 
\begin{eqnarray*}
\lambda_{mn}&=&4\pi^2((m/a)^2+(n/b)^2)\\
&=& 4\pi^2(m^2+\alpha^2n^2) \mbox{ for } (m \in2\mathbb{Z},\, \mathbb{Z}\ni n\geq 0) \mbox{ or } ( m\in2\mathbb{Z}+1,\, \mathbb{Z}\ni n\neq 0)\,. 
\end{eqnarray*}

The number of nodal domains is given by
\begin{equation}\label{nod-3}
\nu_{mn}=2|mn|+\delta_{m0}(|n|+1)+|m|\delta_{n0}\,.
\end{equation}
The first term is due to the checkerboard structure of the nodal set, when $m,n\neq0$. If either $m=0$ or $n=0$ the nodal domains are parallel stripes. This gives rise to the second and third term, where the form of the second term is due to the twisting. Note that $\nu_{mn}$ is defined only for pairs $(m,n)\in\mathbb{Z}^2$, which fulfill the conditions of the definition of the corresponding eigenfunctions above.

\subsection{Flat tori}
The last type of manifolds we will treat are flat tori. We will give the general definition here. Later we will focus on two and three dimensions. For a set of $n$ linearly independent vectors $(\gamma_i)_{1\leq i \leq n}\subset\mathbb{R}^n$ the corresponding lattice is defined by $\Gamma=span_{\mathbb{Z}}\{\gamma_1,\dots,\gamma_n\}$. We obtain a flat torus by $T_{\Gamma}=(\mathbb{R}^n/\Gamma,g_{\Gamma})$, where $g_{\Gamma}=g_0/\Gamma$ and $g_0$ is the standard metric. Like in the previous cases we will reduce the number of free parameters in the class of flat tori in $n$ dimensions by one. In this case we set $\|\gamma_1\|=1$. The resulting class of flat tori in $n$ dimensions is defined by
\[
\mathfrak{Ft}_n=\{T_\Gamma\,: \,\Gamma\mbox{ is an }n\mbox{-dimensional lattice, }\|\gamma_1\|=1\}\,.
\]
The spectrum and the eigenfunctions of the Laplacian may be expressed using the dual lattice $\Gamma^*=\{\gamma^*\in\mathbb{R}^n\,:\,\langle\gamma^*,\gamma\rangle\in\mathbb{Z}\}$ as follows. The eigenvalues of the Laplace operator on a flat torus are given by 
\[
\lambda_{\gamma^*}=4\pi^2\|\gamma^*\|^2,\quad\gamma^*\in\Gamma^*\,.
\]
Explicitly the spectrum is given as follows. For a basis $(g_i^*)_i$ of the dual lattice $\Gamma^*$ of $\Gamma$, given in terms of the standardbasis in $\mathbb{R}^n$ by the representation $g_i^*=(g^1_{i},\dots ,g^n_{i})$, we define the matrix $A=(g^j_i)_{ij}$. The \emph{Gramian matrix}
is then defined by $G=AA^{\top}$. The spectrum of the flat torus $T_{\Gamma}=\mathbb{R}^n\!/\Gamma$ is obtained as follows
\[
\sigma_{T_{\Gamma}}(\triangle)=\{4\pi^2q^{\top}Gq\,|\,q\in\mathbb{Z}^n\}\,.
\]
The eigenspace of an eigenvalue $\lambda\in \sigma_{T_\Gamma}(\triangle)$ is given by
\[
\mathfrak{E}_{\lambda}=span_{\mathbb{C}}\{\mathbf{x} \mapsto \exp(2\pi\langle\gamma^*,\mathbf{x}\rangle)\,|\,\gamma^*\in\Gamma^*,\, 4\pi^2\|\gamma^*\|^2=\lambda\}\,.
\]
We will focus on the real and imaginary parts of the eigenfunctions, that form a basis of these eigenspaces, namely
\[
\varphi^{\Re}_{\gamma^*}(\mathbf{x})=\cos(2\pi\langle\gamma^*,\mathbf{x}\rangle) \mbox{ and } 
\varphi^{\Im}_{\gamma^*}(\mathbf{x})=\sin(2\pi\langle\gamma^*,\mathbf{x}\rangle)\,.
\]
In \cite{BRKL08} a formula for the number of nodal domains of these eigenfunctions was proved. It is identical for both types of eigenfunctions and given by
\begin{equation}
\nu(\gamma^*)=2\mathrm{gcd}(k_1,\hdots,k_n)\,.
\end{equation}
Here $\gamma^*=\sum_ik_i\gamma^*_i$, where $(\gamma^*_i)$ is any basis of the dual lattice $\Gamma^*$ and "$\mathrm{gcd}$" stands for the \emph{greatest common denominator}. It is easy to see that the  number of nodal domains is invariant under a transformation of the basis. 
\subsection{Nodal Sequences}
We will give a definition of the nodal sequences in the specific cases. For every eigenvalue $\lambda_i$ in the ordered spectrum $0=\lambda_0\leq\lambda_1\leq\lambda_2\leq\hdots\leq\lambda_i\leq\hdots$ we define the set of nodal counts of eigenfunctions within the corresponding eigenspace. For rectangles, separable flat tori and Klein bottles we define
\[
\mathcal{N}(\lambda)=\{\nu_{mn}\,:\,\lambda_{mn}=\lambda\}\,.
\] 
For flat tori we define
\[
\mathcal{N}(\lambda)=\{\nu(\gamma^*)\,:\,\gamma^*\in\Gamma^*,\,\lambda_{\gamma^*}=\lambda\}\,.
\]
The nodal sequence is then obtained by ordering these sets according to the growing eigenvalues.
\[
\mathfrak{N}=\{\mathcal{N}(\lambda_1),\mathcal{N}(\lambda_2),\hdots\}
\]


\section{Inverse Problems}

\subsection{Inverse Spectral Problems}
Solving the Inverse Spectral Problem for a Riemannian manifold means deduce its metric from the given spectrum of the Laplacian or show that there are no isospectral, non-isometric pairs within a set of manifolds. This problem was solved for subclasses of surfaces of revolution \cite{BRHE84, ZE98}. For flat tori this problem is solved in two dimensions (cf. \cite{BEGAMA71} p 149 - 151), where an explicit method to reconstruct the dual lattice from the spectrum is given. 
For three-dimensional flat tori it is known that there are no isospectral, non-isometric pairs \cite{SCHI96}. For later use we join the results for flat tori to a first Lemma.
\begin{lem}\label{lem1}
There are no isospectral, non-isometric pairs in $\mathfrak{Ft}_2$ and $\mathfrak{Ft}_3$. A method to reconstruct the metric from the spectrum is known in dimension 2. 
\end{lem}

\subsection{Inverse Nodal Problems}
Solving the Inverse Nodal Problem for a Riemannian manifold means determine its metric by using its nodal sequence or show that there are no pairs of non-isometric manifolds with the same nodal sequence. In every case we assume that the given nodal sequence belongs to a specific type of manifold. We are now able to formulate our main results.
\begin{thm}\label{thm1}
Within the sets $\mathfrak{R}$, $\mathfrak{T}$ and $\mathfrak{K}$ the parameter $\alpha$ is uniquely determined by the nodal sequence. An explicit method to reconstruct $\alpha$ exists. 
\end{thm}
In the next chapter we will prove this theorem by deriving the specific method for every case. In the case of flat tori, the following result can be achieved.
\begin{thm}\label{thm2}
Within the set $\mathfrak{Ft}_2$ a flat torus is explicitly determined by its nodal sequence. A method to explicitly reconstruct the flat torus from its nodal sequence exists. There is no pair of non-isometric flat tori in $\mathfrak{Ft}_3$ with identical nodal sequences. 
\end{thm}
The preceding theorem will be proven with use of Lemma \ref{lem1}. We will reconstruct the spectrum from the nodal sequence and thereby show that they contain similar information.

\section{Proof of the Main Theorem}
The idea of the proof is similar for the different cases. The nodal sequence contains the relative position of the eigenvalue $\lambda_i$ and the set of numbers $\mathcal{N}(\lambda_i)$. While the position allows to deduce information about the relative size of two eigenvalues, the set $\mathcal{N}(\lambda_i)$ gives a restriction on the possible numbers $n,m\in\mathbb{Z}$. Combining these informations carefully it is possible to construct a sequence that converges to the free parameter $\alpha$.


\subsection{Dirichlet Rectangles}
\label{rectangles}
We will prove Theorem \ref{thm1} for the set $\mathfrak{R}$. Let $R_\alpha\in\mathfrak{R}$ and 
\[
\mathfrak{N}(R_\alpha)=\{\mathcal{N}(\lambda_1), \mathcal{N}(\lambda_2),\mathcal{N}(\lambda_3),\dots,\,\mathcal{N}(\lambda_i),\,\dots \}
\]
be the nodal sequence of $R_\alpha$.
By constructing a sequence with limit $\alpha$ using only the nodal sequence we prove the following implication:
\begin{equation}\label{eq1}
\alpha\neq\alpha '\,\Rightarrow\, \mathfrak{N}(R_\alpha)\neq\mathfrak{N}(R_{\alpha '})\,.
\end{equation}
Given a nodal sequence $\mathfrak{N}(R_\alpha)$, we denote the prime numbers by $\mathbb{P}$ and choose $p\in\mathbb{P}\cup\{1\}$. There are exactly two positions $i,j$ with $i<j$ in $\mathfrak{N}(R_\alpha)$, with $p\in\mathcal{N}(\lambda_{i}),\mathcal{N}(\lambda_{j})$. As $\alpha\leq 1$ we conclude
\[
\lambda_i=\lambda_{1p}=1+\alpha^2p^2\leq p^2+\alpha^2=\lambda_{p1}=\lambda_j\,.
\footnote{We omit the prefactor $\pi^2$ of the eigenvalues in the remainder of this note, it has no influence on the results.}
\]
Therefore when we find a nodal count on exactly two positions, we know that the one with the higher index belongs to an eigenvalue of the form $p^2+\alpha^2$, while the other one belongs to an eigenvalue of the form $1+\alpha^2p^2$.
In the general case we have any natural number $N\in\mathbb{N}$ with the prime decomposition
\[
N=p_1\cdot p_2\cdot\hdots\cdot p_k,\, \mbox{ where } p_i\leq p_{i+1},\, p_i\in\mathbb{P}\,.
\]
For any decomposition
\[
N=p_{i_1}\cdot\hdots\cdot p_{i_l}\cdot p_{i_{l+1}}\cdot\hdots\cdot p_{i_{k}}\mbox{ with } 1\leq l\leq k\mbox{ and } \pi(\{1,\hdots,k\})=\{i_1,\hdots,i_k\}, 
\]
where $\pi$ is any permutation, the following relation is valid
\[
N^2+\alpha^2>(p_{i_{1}}\cdot\hdots\cdot p_{i_{l}})^2+\alpha^2(p_{i_{l+1}}\cdot\hdots\cdot p_{i_k})^2.
\]
Therefore we can detect the position of every eigenvalue of the form $N^2+\alpha^2$ for $N\in\mathbb{N}$.
As shown above, we also know the position of every eigenvalue of the form $1+\alpha^2p^2$ for $p\in\mathbb{P}\cup\{1\}$.
With this information we want to determine $\alpha$. We construct two sequences and use the more dense sequence to put bounds on the other one. To every $h\in\mathbb{P}\cup\{1\}$ let $h_{\pm}\in\mathbb{N}$ be defined by
\begin{eqnarray*}
h_+&=&\min\{i\in\mathbb{N} \,:\,1+h^2\alpha^2<i^2+\alpha^2 \}\\
h_-&=&\max\{i\in\mathbb{N} \,:\,i^2+\alpha^2<1+h^2\alpha^2 \}\,.
\end{eqnarray*}
It follows from the definition that
\[
H^h_-:=(h_-^2-1)(h^2-1)^{-1}<\alpha^2<(h_+^2-1)(h^2-1)^{-1}=:H^h_+
\]
for every $h\in\mathbb{P}$.
It suffices to show
\begin{equation}\label{eq3}
H^h_+-H^h_-\overset{h\rightarrow \infty}{\longrightarrow} 0\,.
\end{equation}
Then follows $H^h_\pm\longrightarrow \alpha^2$, which completes the proof. In order to show convergence we decompose
\begin{equation}\label{eq4}
\|H^h_+-H^h_-\|=\left\|(h_+-h_-)(h-1)^{-1}\right\| \left\|(h_++h_-)(h+1)^{-1}\right\|\,.
\end{equation}
We denote by $\lceil s \rceil=\min\{i\in\mathbb{N}\,:\, i\geq s\}$ and $\lfloor s \rfloor=\max\{i\in\mathbb{N}\,:\, i\leq s\}$. Then we can obtain an upper, respectively lower, bound as follows.
\[
h_-\leq h_+\leq h_+':=\lceil h\alpha+1 \rceil
\]
This follows by definition from
\[
(h_+')^2+\alpha^2\geq (h\alpha+1)^2+\alpha^2\geq h^2\alpha^2+1\,.
\]
As lower bound we choose
\[
h_-':=\lfloor h\alpha\rfloor\leq h_-\leq h_+,
\]
because
\[
(h_-')^2+\alpha^2\leq \alpha^2h^2+\alpha^2\leq 1+ \alpha^2h^2\,.
\]
We can now easily deduce the following estimates
\begin{eqnarray*}
\left\|(h_++h_-)(h+1)^{-1}\right\|&\leq&2\left\|h_+(h+1)^{-1}\right\|\\
&\leq&2\left\|\lceil h\alpha+1\rceil(h+1)^{-1}\right\|\\
&\leq&\|(h\alpha+2)(h+1)^{-1}\|\rightarrow 2\alpha
\end{eqnarray*}
and
\begin{eqnarray*}
\left\|(h_+-h_-)(h-1)^{-1}\right\|&\leq&\left\|(h_+'-h_-')(h-1)^{-1}\right\|\\
&=&\left\|(\lceil h\alpha+1 \rceil-\lfloor h\alpha\rfloor)(h-1)^{-1}\right\|\\
&\leq&2(h-1)^{-1}\rightarrow 0\,.
\end{eqnarray*}
With (\ref{eq4}) we obtain (\ref{eq3}) and Theorem \ref{thm1} is proven for $\mathfrak{R}$.$\hfill\Box$


\subsection{Separable Tori}
It is even easier to determine $\alpha$ for separable tori in two dimensions than for the Dirichlet rectangles. This is due to the specific formula for the nodal count (cf. formula  (\ref{nod-2})). Whenever a nodal count of the form $\nu_{mn}=2k$ with $k\in2\mathbb{Z}+1$ appears in the nodal sequence, the eigenvalue in the corresponding position is either $k^2$ or $\alpha^2k^2.$ As $\alpha\leq 1$ it is possible to identify the bigger position as $k^2$ and the smaller as $\alpha^2k^2$. It is therefore possible to determine the position for every eigenvalue of the form $k^2$ or $\alpha^2k^2$ for $k$ odd.
Similar to the previous case we define, for $h\in2\mathbb{Z}+1$,  
\begin{eqnarray*}
h_+&=&\min\{i\in2\mathbb{Z}+1 \,:\,\alpha^2h^2<i^2 \}\\
h_-&=&\max\{i\in2\mathbb{Z}+1 \,:\,i^2<\alpha^2h^2 \}\,.
\end{eqnarray*}
We have
\[
H^h_-:=(h_-/h)^{2}<\alpha^2<(h_+/h)^{2}=:H^h_+\,.
\]
The difference $h_+-h_-$ can be bounded by using
\[
\lfloor h\alpha-2\rfloor<h_-<h_+<\lceil h\alpha+2\rceil\,.
\]
The convergence follows as in the previous case and Theorem \ref{thm1} is proven for $\mathfrak{T}$.$\hfill\Box$


\subsection{Flat Klein bottles}
The Klein bottle case is crucially different to the cases above. This is due to the fact that the equivalence condition is not invariant under a rotation of $\pi/2$. Therefore $K_{\alpha}$ and $K_{\alpha^{-1}}$ are not isometric and $\alpha\in\mathbb{R}_+$. 
To reconstruct $\alpha$ from the nodal sequence it is therefore necessary to distinguish between the cases $\alpha>1$, $\alpha=1$ and $\alpha<1$. In each case it will be possible to apply similar techniques as before. Distinguishing between the cases is possible by analyzing the nodal sequence for low eigenvalues.\\
First of all it is easy to see that whenever an odd nodal count $\nu=(n+1)\in2\mathbb{Z}+1$ occurs, the corresponding eigenvalue has the form $\lambda_{0n}=n^2\alpha^2$ and is simple.\footnote{Here one should recall the restrictions to allowed pairs $(m,n)$ in the Klein bottle case.} We investigate the eigenvalues with nodal count $\nu=2$. Due to (\ref{nod-3}) the following pairs $(m,n)\in\mathbb{Z}$ have nodal count $\nu=2$:
\[
(m_1,n_1)=(1,1),\, (m_2,n_2)=(2,0), \mbox{ and } (m_3,n_3)=(0,1)\,.
\]
 The corresponding eigenvalues are:
\[
 \lambda_{11}=1+\alpha^2,\, \lambda_{20}= 4 \mbox{ and } \lambda_{01}=\alpha^2\,.
 \]
The nodal sequence reveals the positions of this set as a whole, but not the specific positions of its elements. Due to the remark above on odd nodal counts we also know the position of $\lambda_{02}=4\alpha^2$. If the position of $4\alpha^2$ is bigger than the other three, we can conclude $\alpha^2>1$. If the position of $4\alpha^2$ is smaller than exactly one position it follows that
 \[
 4>4\alpha^2>1+\alpha^2>\alpha^2 \mbox{, so } \alpha^2<1\,.
 \] 
 \noindent When $4\alpha^2$ is on the third position, it follows $\alpha^2<1$. The fourth position is impossible. If $4\alpha^2$ shares the maximal position of the set with one element, this element is $4$ and we obtain $\alpha=1$. If it shares the second position, we can conclude $\alpha^2=1/3$. These are all the possibilities, so we can either determine $\alpha$ or decide whether $\alpha^2<(>)1$.\\
\\
We treat the cases $\alpha^2>1$ and $\alpha^2<1$ separately and start with $\alpha^2<1$. We use the following terminology: A relation $A_n\sim B_n$ is \emph{asymptotically valid} (a.v.), iff there is an $n_0\in\mathbb{N}$, so that $\forall n\geq n_0:$ $A_n\sim B_n$ is true.\\
We have seen how the eigenvalues with odd nodal counts  can be expressed in terms of $\alpha^2$. For an even nodal count $\nu=2k$ we are able to identify the positions of the following set of eigenvalues:
\begin{eqnarray*}
&&\lambda_{2k,0}=4k^2,\, \lambda_{0,2k-1}=(2k-1)^2\alpha^2\mbox{ and all values } \lambda_{pq}=p^2+\alpha^2q^2,\\
&&\mbox{where }k=p\cdot q\mbox{ is any decomposition of }k\mbox{ with }p,q\in\mathbb{N}\,.
\end{eqnarray*}
Before $\alpha$ can be reconstructed it is necessary to decide whether $\alpha$ is bigger, smaller or equal to $1/4$. This is important for later use. It will be done by comparing the position of the eigenvalues corresponding to odd nodal counts $\nu=n+1$ - $\lambda_{0n}=n^2\alpha^2$ for $n\in2\mathbb{Z}$ - with a set of eigenvalues for even nodal counts as described above for $\nu=n=2k$. The position of $\lambda_{0n}$ within this set of eigenvalues can be detected by the different nodal counts. It will then be used to deduce information on $\alpha$. We consider $k\in\mathbb{P}$, so the set of eigenvalues with nodal count $\nu=2k$ contains exactly four elements. Due to $\alpha^2<1$ it is clear that $\lambda_{2k,0}>\lambda_{0n}>\lambda_{0,2k-1}$. Furthermore $\lambda_{0n}>\lambda_{1k}$ is asymptotically valid. Whether $\lambda_{0n}\leq\lambda_{k1}$ or $\lambda_{0n}>\lambda_{k1}$ can be read from the nodal sequence, as all other relative positions within the set are known. There are two cases. If $\lambda_{0n}>\lambda_{k1}$ for some $k$, we can conclude $\alpha^2>k^2/(4k^2-1)>1/4$. Otherwise we have $\lambda_{0n}\leq\lambda_{k1}\forall k\in\mathbb{P}$ and therefore $\alpha^2\leq k^2/(4k^2-1)\rightarrow 1/4$.\\
Knowing these bounds we are able to reconstruct $\alpha$ from the nodal sequence.
To identify the specific positions of the elements of the set of eigenvalues with even nodal counts, we need to analyze its structure. As $\alpha^2<1$, the following relations are valid: 
\begin{eqnarray*}
\lambda_{2k,0}=4k^2&>&(2k-1)^2\alpha^2\mbox{ is a.v.}\\
\lambda_{2k-1,0}=(2k-1)^2\alpha^2&>&p^2+\alpha^2q^2 \,\forall p,q :p\cdot q=k  \mbox{ is a.v. if } \alpha>1/4\\ 
\lambda_{2k-1,0}=(2k-1)^2\alpha^2&<& k^2+\alpha^2\mbox{ is a.v. if } \alpha\leq 1/4 \\
\lambda_{k,1}=k^2+\alpha^2&\geq&p^2+\alpha^2q^2  \,\forall p,q :p\cdot q=k\,.
\end{eqnarray*}
It is therefore asymptotically valid that the third (second) largest position among a set of eigenvalues of an even nodal count corresponds to the eigenvalue $k^2+\alpha^2$ if $\alpha^2>1/4$ ($\alpha^2\leq1/4$). For $k$ prime the set consists of four eigenvalues, where the lowest position asymptotically corresponds to the eigenvalue $1+k^2\alpha^2$. We are now able to construct a sequence with limit $\alpha^2$. For any $k\in\mathbb{P}$ we define
\begin{eqnarray*}
k_+&=&\min\{i\in\mathbb{N} \,:\,1+k^2\alpha^2<i^2+\alpha^2 \}\\
k_-&=&\max\{i\in\mathbb{N} \,:\,i^2+\alpha^2<1+k^2\alpha^2 \}\,.
\end{eqnarray*}
We always identify the third (second)\footnote{if $\alpha\leq1/4$} position of a set of eigenvalues of a fixed even nodal count $\nu=2i$ with
$i^2+\alpha^2$. For small nodal counts these positions might be incorrect, but as we construct a sequence for $k\rightarrow\infty$ and the identification is asymptotically valid, this finite set of incorrect values does not affect the limit. The limit is now constructed in the same way as in the rectangle case (cf. page \pageref{rectangles}).\\
\newline
We now study the case $\alpha^2>1$. To know the asymptotic behaviour it is necessary to decide between $\alpha^2>4$ and $\alpha^2<4$. Therefore we investigate the following nodal counts. As mentioned above, we know the position of the eigenvalue $\lambda_{02}=4\alpha^2$. The positions of the set of eigenvalues 
\[
\lambda_{40}=16,\,\lambda_{03}=9\alpha^2,\, \lambda_{12}=4\alpha^2+1\mbox{ and }\lambda_{21}=\alpha^2+4
\] 
are also known. Due to the relation $9\alpha^2>4\alpha^2+1>4\alpha^2$ the position of $4\alpha^2$ might be bigger, smaller or between the positions of $\alpha^2+4$ and $16$. If it is bigger, we can conclude $\alpha^2>4$, if smaller then $\alpha^2<4$. If it lies inbetween, the only possible configuration is $16>4\alpha^2>\alpha^2+4$ and therefore $\alpha^2<4$. If it coincides with the bigger value, we get $\alpha^2=4$, if it coincides with the smaller we get $\alpha^2=4/3$. In any case, we are able to determine $\alpha^2$ or decide between $\alpha^2>4$ and $\alpha^2<4$. We have the following relations, due to $\alpha^2>1$:
\begin{eqnarray*}
\lambda_{0,2k-1}=(2k-1)^2\alpha^2&>&4k^2\mbox{ is a.v.}\\
\lambda_{0,2k-1}=(2k-1)^2\alpha^2&>&p^2+\alpha^2q^2\,\forall p,q : p\cdot q=k \mbox{ is a.v.}\\
\lambda_{1,k}=1+k^2\alpha^2&\geq&p^2+\alpha^2q^2\,\forall p,q : p\cdot q=k\\
\lambda_{2k,0}=4k^2&>&\alpha^2+k^2\mbox{ is a.v.}
\end{eqnarray*}
As above, we can therefore asymptotically determine the positions of the eigenvalues of the form $\alpha^2+k^2$ for $k\in\mathbb{P}$ as the smallest position in a given set of eigenvalues with the same nodal count $\nu=2k$. Furthermore we have another asymptotically valid relation depending on $\alpha^2$:
\begin{eqnarray*}
4k^2>1+\alpha^2k^2&&\mbox{ is a.v. if }4>\alpha^2\\
4k^2<1+\alpha^2k^2&&\mbox{ is a.v. if }4<\alpha^2\,.  
\end{eqnarray*}
With this we can asymptotically determine the position of eigenvalues of the form $1+\alpha^2k^2$ as the third (second) biggest position in a set of eigenvalues with nodal count $\nu=2k$ for any $k$ if $\alpha^2>4$ ($\alpha^2<4$). In this case ($\alpha^2>1$) the situation is different, since we can detect the positions of $1+\alpha^2k^2$ for any $k$, but the positions $k^2+\alpha^2$ only for $k\in\mathbb{P}$. To construct the sequence we need to change the previous method slightly. For $k$ prime we define
\begin{eqnarray*}
k_+&=&\min\{i\in\mathbb{N} \,:\,k^2+\alpha^2<1+i^2\alpha^2 \}\\
k_-&=&\max\{i\in\mathbb{N} \,:\,1+i^2\alpha^2<k^2+\alpha^2 \}\,.
\end{eqnarray*}
Then for every $k\in\mathbb{P}$ we have the upper and lower bounds
\[
1+\alpha^2k_-^2<\alpha^2+k^2<1+\alpha^2k_+^2,
\] 
which leads to
\[
(k_-^2-1)(k^2-1)^{-1}<\alpha^{-2}<(k_+^2-1)(k^2-1)^{-1}.
\]
After showing, as before, that both sequences converge, we can determine $\alpha^2$ and thus Theorem \ref{thm1} is proved for $\mathfrak{K}$.$\hfill\Box$


\subsection{Flat Tori}

Every eigenvalue $\lambda$ is determined by its representation vectors $(\gamma^*_j)$. For a given basis $(\gamma_i)_{i\leq n}$ of $\Gamma$ and the corresponding dual basis each representation vector is given by a coefficient vector $\mathrm{q}_j\in\mathbb{Z}^n$, where $n$ is the dimension of the flat torus. The greatest common denominator of the coefficient vector is invariant under transformations of the basis. That means we can associate this $\mathrm{gcd}$ directly with the representation vector and call it the \emph{$\mathit{gcd}$ of $\gamma^*$}. For each representation vector of an eigenvalue there is one corresponding gcd. Therefore the number of greatest common denominators of an eigenvalue equals its degeneracy.\\
\noindent The complete spectrum is already determined by all eigenvalues, which have at least one 
greatest common denominator that equals one. All other eigenvalues are obtained by multiplying these eigenvalues with squares of natural numbers.
\[
\lambda=4\pi^2q\!^{\top}\!Qq=\mathrm{gcd}(q)^24\pi^2\widetilde{q}^{\top}\!Q\widetilde{q}=\mathrm{gcd}(q)^2\widetilde{\lambda}\mbox{ with }\mathrm{gcd}(\widetilde{q})=1
\]
We call the eigenvalues $\widetilde{\lambda}$ \emph{basic eigenvalues}. Knowing the nodal sequence we are able to find the positions of all basic eigenvalues in the ordered spectrum, by choosing the positions which have nodal count $\nu=2$. The set of these positions is $\mathcal{P}_2=\{m\in\mathbb{N}\, |\,2\in\mathcal{N}(\lambda_{m})\}$.
In the same way we define the set of positions of the $k^2$-multiples for $k\in\mathbb{N}$ of the basic eigenvalues 
\[
\mathcal{P}_k=\{m\in\mathbb{N}\, |\,2k\in\mathcal{N}(\lambda_m)\}\,.
\]
We denote by $p_k^i$ the $i$-th element of $\mathcal{P}_k$. The set of positions of a sequence of eigenvalues which are the $k^2$-multiples of the $i$-th basic eigenvalue is given by $(p_k^i)_{k\in\mathbb{N}}$.
That is because multiplying two eigenvalues with a square number conserves their order. Therefore, the nodal count $2k$ that appears at $m\in\mathbb{N}$ at the $i$-th position in $\mathcal{P}_k$ indicates that the eigenvalue $\lambda_m$ equals $k^2\widetilde{\lambda}_i$, where $\widetilde{\lambda}_i$ is the $i$-th basic eigenvalue in the ordered spectrum. 
Using the relative positions of the eigenvalues we want to reconstruct their ratio. Therefore we choose two arbitrary sequences $(p^i_k)_{k\in\mathbb{N}}$ and $(p^j_l)_{l\in\mathbb{N}}$. These two sequences represent the relative positions of the basic eigenvalues $\widetilde{\lambda}_i$ and $\widetilde{\lambda_j}$ and their multiples with squares of natural numbers. We will construct a sequence whose limit equals $\widetilde{\lambda_j}/\widetilde{\lambda_i}$.\\
\noindent 
Without loss of generality let be $p^i_1\leq p^j_1$. Since $p^i_1=p^j_1$ implies $\lambda_j/\lambda_i=1$, we consider the case $p^i_1< p^j_1$ from now on. That means $p^i_m< p^j_m$ $\forall m\in\mathbb{N}$. For a given $m\in\mathbb{N}$ we define $m_{\max}=\max\{k\in\mathbb{N}\,|\,p^i_k\leq p^j_m\}$. From the definition it follows $p^i_{m_{\max}}\leq p^j_{m}\leq p^i_{m_{\max}+1}$.
In terms of the eigenvalues this becomes
\[
m_{\max}^2\widetilde{\lambda}_i\leq m^2\widetilde{\lambda_j}\leq (m_{\max}+1)^2\widetilde{\lambda}_i\,.
\] 
We obtain the following bounds for the ratio:
\[
m_{\max}^2/m^2\leq\widetilde{\lambda}_j/\widetilde{\lambda}_i\leq (m_{\max}+1)^2/m^2\,.
\]
As $m\rightarrow\infty$, $m_{\max}\rightarrow\infty$ and we get 
\begin{eqnarray*}
\|(m_{\max})^2/m^2-(m_{\max}+1)^2/m^2\|&\leq&2\|m_{\max}/m^2|+|1/m^2\|\\
&=&2\underbrace{\|(m_{\max})^2/m^2\|}_{\leq\widetilde{\lambda}_j/\widetilde{\lambda}_i}\|1/m_{\max}\|+\|1/m^2\|\overset{m\rightarrow\infty}{\longrightarrow}0\,.
\end{eqnarray*}
We can now deduce
\[
\lim_{m\rightarrow\infty}(m_{\max})^2/m^2=\widetilde{\lambda}_j/\widetilde{\lambda}_i\,.
\]
In this way we can determine the size of all basic eigenvalues up to a factor. The complete spectrum
is given by the basic eigenvalues multiplied with squares of natural numbers. Applying Lemma \ref{lem1} completes the proof of Theorem \ref{thm2}. $\hfill\Box$ 

\section{Remarks}
The examples that we have studied support the conjecture on the role of nodal sequences in the determination of manifolds. The principal approach to determine a manifold by its nodal sequence is similar for all classes presented here. Even though the investigated systems all have a very simple structure, it might  be interesting to see whether the main idea - identifying positions in the spectrum with explicit eigenvalues by using the corresponding nodal counts -  can be generalized for more complex systems. Using similar ideas for the case of surfaces of revolution that are treated in \cite{karasmil08} it is possible to deduce some basic properties of the surface \cite{SMpr}.
\section*{Acknowledgements}
The author thanks Rami Band for enlightening discussions on the area of nodal domains, some of which led to this work, and Amit Aronovitch for several conversations. The author is grateful to Uzy Smilansky for his strong encouragement, support and suggestions. The author also wants to thank the Weizmann Institute of Sciences for the hospitality during his visits, and the participants to the Quantum Chaos seminar at the Institute of Physics of Complex Systems at the Weizmann Institute, especially Yehonatan Elon, Idan Oren and Nir Auerbach. The author acknowledges the financial support of the German Israeli Foundation (GIF).


\end{document}